\documentclass{article}%[article]%{smfclass}

\usepackage[utf8]{inputenc}
\usepackage[T1]{fontenc}
\usepackage{amsmath, amssymb, amsthm}
\usepackage{url}
\usepackage{hyperref}
\usepackage{graphicx}
\usepackage{natbib}
\usepackage{lmodern}   % ou libertinus
\usepackage{microtype}

\usepackage{tikz}
\usetikzlibrary{positioning}

\newcommand{\permutationset}[1]{\mathcal{S}_{#1}}          % permutations sur [#1]
\newcommand{\arrangementset}[2]{\mathcal{A}_{#1}^{(#2)}}   % injections [#1] -> [#2]
\newcommand{\graphset}[2]{\mathcal{G}_{#1,#2}}             % \ell-appariements (n, ell)

\newcommand{\derangements}[1]{D_{#1}}            % D_n : dérangements (permutations sans point fixe)
\newcommand{\rectderangements}[2]{D_{#1}^{(#2)}} % D_n^{(m)} : dérangements rectangulaires (injections sans point fixe)
\newcommand{\partialderangements}[2]{D_{#1,#2}}  % D_{n,\ell} : dérangements partiels (\ell-appariements sans point fixe)

\newcommand{\grapharrangementsset}[3]{\mathcal{G}_{#1,#3}^{(#2)}}
\newcommand{\rectpartialderangements}[3]{D_{#1,#3}^{(#2)}}
\newcommand{\Prectpartial}[3]{P_{#1,#3}^{(#2)}}  % proba associée (optionnel)

\newenvironment{keywords}{%
  \par\vspace{.5\baselineskip}\noindent\small
  \textbf{Keywords}—\ignorespaces
}{\par}

\title{Derangements and Generalizations: A Counting Note on the Matching Problem}
\author{Antoine Luciano}
\date{\today}

\begin{document}

\maketitle
  
\begin{abstract}
We give a concise historical background to \emph{Montmort’s matching problem} and its modern variants such as the \emph{hat-check problem}, then develop a unified counting framework for fixed-point-free allocations. Using elementary recurrence and inclusion–exclusion arguments, we derive closed forms for derangements, rectangular injections, and partial $\ell$-matchings, and we combine them into a single formula. We also provide exact counts for the number of fixed points and Poisson limit laws. This note thus offers a compact, self-contained synthesis linking classical results with their two principal generalizations in a single scheme.
\end{abstract}

\begin{keywords}
derangements; matching problem; rencontres numbers; partial matchings; \(\ell\)-matchings; inclusion--exclusion principle; recurrence relations; permutations; injections
\end{keywords}

\begin{figure}
    \centering
    \includegraphics[width=0.9\textwidth]{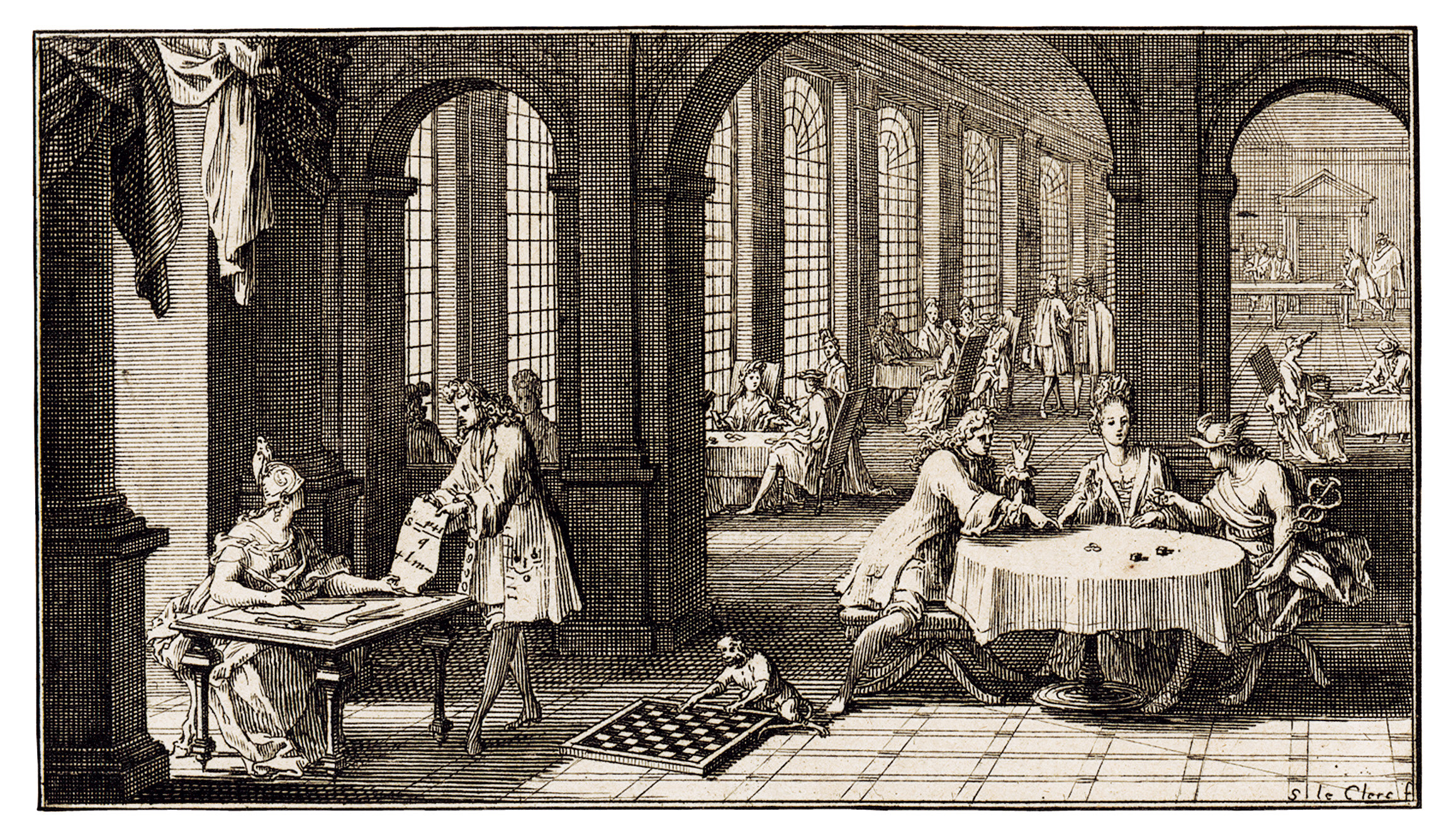} 
    \caption{Engraving by Sébastien Leclerc from the Preface of Pierre Rémond de Montmort's \textit{Essay d'Analyse sur les Jeux de Hazard} \citet{demontmort1708essay}.}
    \label{fig:montmort}
\end{figure}

\section*{Introduction}

The \emph{matching problem}, also known as \emph{Montmort's problem} or the \emph{hat-check problem}, is a classic of combinatorics and probability theory. Formulated in various forms over more than three centuries, it provides an elegant thread connecting card games, urn drawings, chessboard placements, and permutations.

Its origin goes back to the first edition of the \textit{Essay d’Analyse sur les Jeux de Hazard} by \citet{demontmort1708essay}, which presents a simplified version of the game of \textit{Treize} (\emph{Thrirteen}): a player calls out the values from 1 to 13 (Ace, 2, \dots, King) while drawing cards from a full deck; the player wins if, at the moment of the call, the card of the corresponding value is drawn, and loses if no match occurs after all thirteen draws. The solution to the problem was provided by \citet{demontmort1713essay}, in correspondence with Nicolas Bernoulli included in the second edition of the book.

In his \textit{Calcul de la probabilité dans le jeu de rencontre}, \citet{euler1753calcul} studied a variant in which two players each draw a card from a deck of \(n\) cards numbered from 1 to \(n\). If the two cards have the same number, the first player wins; otherwise, the second does. Euler then computed the probability that the two cards differ.

\citet{catalan1837solution} reformulated the problem in terms of urn drawings with numbered balls and generalized it: determining the probability that exactly \(k\) objects (in this setting, balls) return to their initial position.

In his \textit{Théorie des nombres}, \citet{lucas1891theorie} formulated an analogue known as the \textit{rook problem}: placing \(n\) rooks on an \(n \times n\) chessboard with exactly one rook per row and per column, with none lying on the main diagonal.

The most popular modern version is the \textit{hat-check problem}, popularized by \citet{scoville1966hat}: \(n\) guests leave their hats in a cloakroom, which are then randomly redistributed; the question is the probability that no guest retrieves their own hat.

All these formulations point to the same object: the probability that a random permutation of \(n\) elements is a \emph{derangement}, that is, has no fixed point. Despite its apparent simplicity, the problem admits numerous generalizations and applications and continues to be studied today \citep{capparelli2018generalization}. It has also emerged in the context of my doctoral research in Bayesian statistics \citep{luciano2025permutations}. As I have thus far not found a survey note gathering these perspectives, I propose here a unified presentation of the basic results and their extensions, in the hope that it may be useful to others.

% \begin{figure}[ht]
%     \centering
%     \includegraphics[width=0.9\textwidth]{scoville.png} 
%     \caption{Illustration from \textit{The Hat Check Problem} \citep{scoville1966hat}}
%     \label{fig:scoville}
% \end{figure}

\section{Initial Problem: Permutations and Derangements}
\label{sec:perm}

Throughout this paper, we will use the \textit{hat-check problem} framework described earlier \citep{scoville1966hat}. To formalize the problem mathematically, let \(\permutationset{n}\) denote the set of permutations of \(\{1,2,\dots,n\}\), that is, bijections from \(\{1,2,\dots,n\}\) onto itself. We have \(\lvert \permutationset{n}\rvert=n!\), where \(\lvert A \rvert\) denotes the \emph{cardinality} (number of elements) of the set \(A\).

Each \(\sigma\in\permutationset{n}\) represents a redistribution of the hats among the \(n\) people. We are interested in the number of permutations without fixed points (no person receives their own hat):
\[
\derangements{n} \;=\; \bigl\{\,\sigma\in\permutationset{n}:\ \forall i\in\{1,\dots,n\},\ \sigma(i)\neq i\,\bigr\}\,\text{.}
\]

The probability that a uniform random permutation has no fixed points is thus
\[
P_n=\mathbb{P}(\sigma \in \derangements{n})=\frac{\lvert \derangements{n}\rvert}{n!}\,\text{.}
\]

The number \(\lvert \derangements{n}\rvert\) is the \emph{number of derangements} (or \emph{subfactorial}) of \(n\), often denoted \(!n\) or \(d_n\) \citep{OEIS:A000166}. This quantity, first discovered in the course of solving the problem, admits the explicit formula
\begin{equation*}\label{eq:derangement-formula}
\lvert \derangements{n}\rvert \;=\; !n \;=\; n!\sum_{k=0}^{n}\frac{(-1)^k}{k!}
\;=\;  \left\lfloor \frac{n!}{e}\right \rceil \,\text{,}
\end{equation*}
where $\lfloor x \rceil$ denotes rounding $x$ to the nearest integer. Hence
\begin{equation}\label{eq:derangement-proba}
P_n=\sum_{k=0}^{n}\frac{(-1)^k}{k!}\,\text{,}
\qquad\text{and}\qquad
\lim_{n\to\infty}P_n=\frac1e\simeq 0.3679\,\text{.}
\end{equation}
This asymptotic behavior was already observed by \citet{demontmort1713essay}, although the constant \(e\) was only formalized later by \citet{euler1736mechanica}.

\subsection{Recursive Approach}

The recursive method is without doubt the most natural one: we follow the first move and unfold the consequences. This is in fact the reasoning by which the original problem was historically solved by \citet{demontmort1713essay} (based on the notes of Nicolas Bernoulli).

Let us number the people and hats from \(1\) to \(n\).
Consider what person \(1\) does in a permutation without fixed points. She cannot take hat \(1\), so she must take some hat \(j\neq 1\). There are \(n-1\) possible choices for \(j\).

Two situations then arise:

\begin{itemize}
  \item \textbf{If person \(j\) takes hat \(1\).}\\
  Persons \(1\) and \(j\) are “paired” (a short cycle of length 2).
  We can then ignore these two people and their hats: the problem reduces
  to \(n-2\) people. This gives \( \lvert \derangements{n-2}\rvert \,\text{.}\)

  \item \textbf{If person \(j\) does not take hat \(1\).}\\
  In this case, we can reduce the situation to a problem with \(n-1\) people
  (removing person \(1\) and hat \(j\), and reindexing the assignments).
  We then recover exactly a derangement of size \(n-1\), yielding \( \lvert \derangements{n-1}\rvert \,\text{.}\)
\end{itemize}

Since there are \(n-1\) choices for \(j\), we obtain the recurrence
\[
\lvert \derangements{n}\rvert=(n-1)\bigl(\lvert \derangements{n-1}\rvert+\lvert \derangements{n-2}\rvert\bigr)\,\text{,}
\qquad \lvert \derangements{0}\rvert=1\,\text{,}\quad \lvert \derangements{1}\rvert=0\,\text{.}
\]

We can also write
\[
\lvert \derangements{n}\rvert=n\,\lvert \derangements{n-1}\rvert+(-1)^n\,\text{,}
\]
and, by dividing by \(n!\) and summing, we recover
\[
P_n=\sum_{k=0}^{n}\frac{(-1)^k}{k!}\,\text{,}
\quad\text{and}\quad
\lvert \derangements{n}\rvert=n!\sum_{k=0}^{n}\frac{(-1)^k}{k!}\,\text{.}
\]

\subsection{Inclusion–Exclusion Approach}

We use the \emph{principle of inclusion–exclusion} (also known as the \emph{Poincaré sieve}). 
For finite sets \(A_1,\dots,A_n\),
\begin{equation}\label{eq:crible-poincare}
\left\lvert \bigcup_{i=1}^{n} A_i \right\rvert
\;=\!\!\sum_{\varnothing\neq S\subseteq [n]}\!\!(-1)^{\lvert S\rvert+1}\left\lvert \bigcap_{i\in S} A_i \right\rvert
\;=\sum_{k=1}^{n}(-1)^{k+1}\!\!\!\sum_{1\le i_1<\cdots<i_k\le n}\!\!\!
\left\lvert A_{i_1}\cap\cdots\cap A_{i_k}\right\rvert \,\text{.}
\end{equation}

Let us apply \eqref{eq:crible-poincare} to the hat-check problem. 
For each \(i\in\{1,\dots,n\}\), set
\[
A_i=\{\sigma\in\permutationset{n}:\ \sigma(i)=i\}\,\text{,}
\]
the set of permutations where person \(i\) retrieves their own hat.
If \(S\subseteq\{1,\dots,n\}\) with \(\lvert S\rvert=k\), then fixing \(\sigma(i)=i\) for all \(i\in S\) leaves \((n-k)!\) free permutations of the remaining indices:
\[
\left\lvert \bigcap_{i\in S}A_i \right\rvert=(n-k)!\,\text{.}
\]
Since there are \(\binom{n}{k}\) subsets \(S\) of size \(k\), we obtain
\[
\left\lvert \bigcup_{i=1}^{n} A_i \right\rvert
=\sum_{k=1}^{n}(-1)^{k+1}\binom{n}{k}(n-k)!
= n!\sum_{k=1}^{n}\frac{(-1)^{k+1}}{k!}\,\text{.}
\]

The derangements are the complement of this union:
\[
\derangements{n}=\permutationset{n}\setminus \bigcup_{i=1}^{n} A_i\,\text{,}
\qquad
\lvert \derangements{n}\rvert
=\lvert \permutationset{n}\rvert-\left\lvert \bigcup_{i=1}^{n} A_i \right\rvert
= n!\sum_{k=0}^{n}\frac{(-1)^k}{k!}\,\text{,}
\]
where the empty term \(k=0\) (equal to \(1\)) is included. Dividing by \(n!\), we recover
\[
P_n=\sum_{k=0}^{n}\frac{(-1)^k}{k!}\,\text{.}
\]

\subsection{Distribution of the Number of Fixed Points (Encounter Numbers)}

We are interested in the probability that \emph{exactly \(k\)} people retrieve their own hat (i.e., \(k\) fixed points).
One must first choose these \(k\) people (\(\binom{n}{k}\) ways), then derange the remaining \(n-k\) (\(\lvert \derangements{n-k}\rvert\) ways).
The total number of permutations of size \(n\) with exactly \(k\) fixed points, also called the \textit{encounter number} or the number of \textit{partial derangements} \citep{OEIS:A047920}, is therefore
\[
N_n(k)=\binom{n}{k}\,\lvert \derangements{n-k}\rvert\,\text{,}\qquad
k=0,1,\dots,n\,\text{.}
\]

We note that \(N_n(n)=1\) (only the identity permutation has \(n\) fixed points), \(N_n(n-1)=0\) (if \(n-1\) people get their own hat, the last one must as well), and \(N_n(0)=\lvert \derangements{n}\rvert\,\text{.}\)

The corresponding probability distribution of the number of fixed points \(K\) (under the uniform distribution on \(\permutationset{n}\)) is, for any fixed \(k\in\{0,\dots,n\}\),
\[
P_n(k)=\frac{N_n(k)}{n!}
=\frac{\binom{n}{k}\,\lvert \derangements{n-k}\rvert}{n!}
=\frac{1}{k!}\sum_{j=0}^{n-k}\frac{(-1)^j}{j!}
\;\approx\; \frac{e^{-1}}{k!}\quad(n\ \text{large})\,\text{.}
\]

In other words, \(K \leadsto \mathrm{Poisson}(1)\). The Poisson approximation can be justified by the representation \(K=\sum_{i=1}^{n}\mathbf 1_{\{\sigma(i)=i\}}\) together with Chen–Stein type arguments.

Finally, the decomposition by the number of fixed points gives
\[
\sum_{k=0}^{n} N_n(k)=n!\,\text{,}\qquad\text{hence } \sum_{k=0}^{n} P_n(k)=1\,\text{.}
\]
\section{First Generalization: Rectangular Case — Fixed-Point-Free Injections}

Suppose now that the cloakroom contains \(m>n\) hats (with $m-n$ “extra” hats without owners).  
We seek the probability that none of the $n$ people retrieves their own hat, which generalizes the case \(m=n\).
This framework appears under the name of \emph{\((n,m)\)-matching} or \emph{fixed-point-free arrangements} \citep{hanson1983matchings}.

We consider the set of \emph{arrangements} of \(n\) objects among \(m\), that is, injections \(\{1,\dots,n\}\hookrightarrow\{1,\dots,m\}\), denoted \(\arrangementset{n}{m}\).
The total number of arrangements is
\[
\lvert \arrangementset{n}{m}\rvert=\frac{m!}{(m-n)!}\,\text{.}
\]

We then define the set of \emph{rectangular derangements} (fixed-point-free injections):
\[
\rectderangements{n}{m}=\Bigl\{\sigma\in\arrangementset{n}{m}:\ \forall i\in\{1,\dots,n\},\ \sigma(i)\neq i\Bigr\}\,\text{.}
\]
By inclusion–exclusion, we obtain the explicit formula \citep{OEIS:A076731}
\begin{equation}\label{eq:rect-derangement-formula}
\lvert \rectderangements{n}{m}\rvert=\sum_{k=0}^{n}(-1)^k\binom{n}{k}\,\frac{(m-k)!}{(m-n)!}\,\text{.}
\end{equation}

Thus, for a uniformly chosen injection in \(\arrangementset{n}{m}\),
\[
P_n^{m}
=\frac{\lvert \rectderangements{n}{m}\rvert}{\lvert \arrangementset{n}{m}\rvert}
=\sum_{k=0}^{n}(-1)^k\binom{n}{k}\,\frac{(m-k)!}{m!}\,\text{.}
\]

Note that when \(m=n\), we recover the formulas of the initial case (Section~\ref{sec:perm}).

\subsection{Inclusion–Exclusion Approach (Poincaré Sieve)}

We adapt the Poincaré sieve to the rectangular case.
For \(i\in\{1,\dots,n\}\), set
\[
A_i=\{\sigma\in\arrangementset{n}{m}:\ \sigma(i)=i\}\,\text{.}
\]
Fixing simultaneously \(k\) indices (that is, imposing \(\sigma(i)=i\) for \(i\) in a subset \(S\) of size \(k\)) leaves \(n-k\) elements of the domain to be injected into a codomain of size \(m-k\), hence
\[
\Bigl\lvert\bigcap_{i\in S}A_i\Bigr\rvert=\lvert \arrangementset{n-k}{m-k}\rvert=\frac{(m-k)!}{(m-n)!}\,\text{.}
\]
Applying \eqref{eq:crible-poincare}, we obtain
\[
\left\lvert \bigcup_{i=1}^{n}A_i\right\rvert
=\sum_{k=1}^{n}(-1)^{k+1}\binom{n}{k}\,\frac{(m-k)!}{(m-n)!}\,\text{,}
\]
and therefore
\[
\lvert \rectderangements{n}{m}\rvert=\lvert \arrangementset{n}{m}\rvert-\left\lvert \bigcup_{i=1}^{n}A_i\right\rvert
=\sum_{k=0}^{n}(-1)^k\binom{n}{k}\,\frac{(m-k)!}{(m-n)!}\,\text{.}
\]

\subsection{Distribution of the Number of Fixed Points (Encounter Numbers)}

We now seek the number of injections with \emph{exactly \(k\)} fixed points.
We first choose these \(k\) indices (\(\binom{n}{k}\) ways), then impose that none of the remaining \(n-k\) indices is fixed in the reduced codomain of size \(m-k\).
In other words, we count \emph{rectangular derangements} of size \((n-k,m-k)\).
We obtain
\[
N_{n}^{m}(k)=\binom{n}{k}\,\lvert \rectderangements{n-k}{m-k}\rvert\,\text{,}
\qquad
k=0,1,\dots,n\,\text{.}
\]

We note that 
\[
N_{n}^{m}(n)=1\quad(\text{the identity injection})\,\text{,} 
\qquad
N_{n}^{m}(n-1)=n\,(m-n)\,\text{,}
\]
(since one index remains to be sent outside itself, with \(m-n\) choices),
and \(N_{n}^{m}(0)=\lvert \rectderangements{n}{m}\rvert\,\text{.}\)

The corresponding probability distribution of the number of fixed points \(K\) (under the uniform distribution on \(\arrangementset{n}{m}\)) is, for any fixed \(k\in\{0,\dots,n\}\),
\[
P_{n}^{m}(k)=\frac{N_{n}^{m}(k)}{\lvert \arrangementset{n}{m}\rvert}
= \frac{n!}{m!}\,\frac{1}{k!}\sum_{j=0}^{n-k}\frac{(-1)^j}{j!}\,\frac{(m-k-j)!}{(n-k-j)!}\,\text{.}
\]

Thus, when \(n,m\to\infty\) with \(n/m\to\rho\in(0,1)\), the number of fixed points \(K\) converges in distribution to \(\mathrm{Poisson}(\rho)\).  
In particular,
\[
P_{n}^{m}(k)= \mathbb P(K=k)\ \xrightarrow{n,m\to\infty}\ e^{-\rho}\,\frac{\rho^{k}}{k!}\qquad(k\ \text{fixed})\,\text{.}
\]

Finally, the decomposition by the number of fixed points gives
\[
\sum_{k=0}^{n}N_{n}^{m}(k)=\lvert \arrangementset{n}{m}\rvert=\frac{m!}{(m-n)!}\,\text{,}
\qquad
\sum_{k=0}^{n}P_{n}^{m}(k)=1\,\text{.}
\]

\section{Second Generalization: Partial Matchings — Fixed-Point-Free $\ell$-Matchings}

We now consider the case where only \(\ell\le n\) people leave with a hat among the \(n\) available.  
This framework is studied under the name of \emph{\(\ell\)-matchings} or \emph{partial injections} \citep{hanson1983matchings}.
An \emph{\(\ell\)-matching} is an injective assignment of \(\ell\) distinct guests to \(\ell\) distinct hats, with no constraint on the remaining \(n-\ell\) guests (who do not receive a hat).
Formally, we consider the set of \emph{partial injections} of \(\{1,\dots,n\}\) into itself, denoted \(\graphset{n}{\ell}\):
\[
\graphset{n}{\ell}
=\bigl\{\,\sigma:I\hookrightarrow \{1,\dots,n\}\ \text{injection}\ :\ I\subseteq \{1,\dots,n\},\ |I|=\ell\,\bigr\}\,\text{.}
\]

To compute \(\lvert \graphset{n}{\ell}\rvert\), proceed in three steps:  
(i) choose the \(\ell\) people among \(n\) (\(\binom{n}{\ell}\) ways),  
(ii) choose the \(\ell\) hats among \(n\) (\(\binom{n}{\ell}\) ways),  
(iii) choose a bijection between these two subsets (\(\ell!\) ways).  
Thus
\[
\lvert \graphset{n}{\ell}\rvert=\binom{n}{\ell}^{2}\,\ell!\,\text{.}
\]

One way to visualize this problem is through the \emph{complete bipartite graph} \(K_{n,n}\):  
the \(n\) people are the vertices on the left, the \(n\) hats those on the right, and an edge \((i,j)\)  
indicates that person \(i\) receives hat \(j\).  
An \(\ell\)-matching is then a set of \(\ell\) pairwise disjoint edges (at most one edge incident to each vertex).  
\emph{Fixed points} correspond to “diagonal” edges \((i,i)\).  
In this language, one speaks of \(\ell\)-matchings with exactly \(k\) fixed points (exactly \(k\) diagonal edges),  
or \emph{fixed-point-free} (no diagonal edge).

We thus define the set of \emph{fixed-point-free \(\ell\)-matchings}:
\[
\partialderangements{n}{\ell}
=\bigl\{\,\sigma\in\graphset{n}{\ell}:\ \forall i\in I,\ \sigma(i)\neq i\,\bigr\}\,\text{.}
\]
By inclusion–exclusion (Poincaré sieve), we obtain \citep{OEIS:A002467}
\begin{equation*}\label{eq:partial-derangements}
\lvert \partialderangements{n}{\ell}\rvert
=\sum_{j=0}^{\ell}(-1)^j\binom{n}{j}\,\binom{n-j}{\ell-j}^{2}\,(\ell-j)!\,\text{.}
\end{equation*}

Under the uniform distribution on \(\graphset{n}{\ell}\),
\[
P_{n,\ell}
=\mathbb{P}(\sigma\in \partialderangements{n}{\ell})
=\frac{\lvert \partialderangements{n}{\ell}\rvert}{\lvert \graphset{n}{\ell}\rvert}
=\frac{1}{\binom{n}{\ell}^{2}\,\ell!}
\sum_{j=0}^{\ell}(-1)^j\binom{n}{j}\,\binom{n-j}{\ell-j}^{2}\,(\ell-j)!\,\text{.}
\]

Note that when \(\ell=n\), we have \(\lvert \graphset{n}{n}\rvert=n!\) and \(\lvert \partialderangements{n}{n}\rvert={!n}\),  
which recovers the initial case of Section~\ref{sec:perm}.

\subsection{Inclusion–Exclusion Approach (Poincaré Sieve)}

For \(i\in\{1,\dots,n\}\), set \(A_i=\{\sigma\in\graphset{n}{\ell}:\ \sigma(i)=i\}\,\text{.}\)  
If we force \(k\) indices to be fixed (a subset \(S\) of size \(k\)), it remains to choose \(\ell-k\) people among \(n-k\), \(\ell-k\) hats among \(n-k\), then a bijection, giving
\[
\Bigl\lvert\bigcap_{i\in S}A_i \Bigr\rvert
=\binom{n-k}{\ell-k}^{2}\,(\ell-k)!\,\text{.}
\]
Applying the inclusion–exclusion principle from Equation~\eqref{eq:crible-poincare} with alternating signs recovers the result of Equation~\eqref{eq:partial-derangements}\,\text{.}

\subsection{Distribution of the Number of Fixed Points (Encounter Numbers)}
f
We now study the number of partial \(\ell\)-matchings with \textit{exactly \(k\)} fixed points, denoted \(N_{n,\ell}(k)\).  
Once again, we first choose the \(k\) fixed indices (\(\binom{n}{k}\) ways), then forbid any further fixed point among the remaining \(\ell-k\) matched indices, which corresponds to an \((\ell-k)\)-matching of \(n-k\) hats without fixed points.  
We obtain
\[
N_{n,\ell}(k)=\binom{n}{k}\,\lvert \partialderangements{n-k}{\ell-k}\rvert\,\text{,} \qquad
k=0,1,\dots,\ell\,\text{.}
\]

In particular,
\[
N_{n,\ell}(0)=\lvert \partialderangements{n}{\ell}\rvert\,\text{,}\qquad
N_{n,\ell}(\ell)=\binom{n}{\ell}\,\text{,}\qquad
N_{n,\ell}(\ell-1)=\binom{n}{\ell-1}\,(n-\ell+1)(n-\ell)\,\text{.}
\]

We deduce the corresponding probability distribution of the number of fixed points \(K\) (under the uniform distribution on \(\graphset{n}{\ell}\)), for any fixed \(k\in\{0,\dots,\ell\}\),
\[
P_{n,\ell}(k)=\frac{N_{n,\ell}(k)}{\lvert \graphset{n}{\ell}\rvert}
= \frac{\ell!}{n!}\,\frac{1}{k!}\sum_{j=0}^{\ell-k}\frac{(-1)^j}{j!}\,\frac{(n-k-j)!}{(\ell-k-j)!}\,\text{.}
\]

Thus, when \(n\to\infty\) with \(\ell/n\to\theta\in[0,1]\), the number of fixed points \(K\) of a uniform \(\ell\)-matching converges in distribution to \(\mathrm{Poisson}(\theta)\).  
In particular,
\[
P_{n,\ell}(k)= \mathbb P(K=k)\ \xrightarrow{n\to\infty}\ e^{-\theta}\,\frac{\theta^{k}}{k!}\qquad(k\ \text{fixed})\,\text{.}
\]

Finally, the decomposition by the number of fixed points gives
\[
\sum_{k=0}^{\ell}N_{n,\ell}(k)=\lvert \graphset{n}{\ell}\rvert
=\binom{n}{\ell}^{2}\,\ell!\,\text{,}
\qquad
\sum_{k=0}^{\ell}P_{n,\ell}(k)=1\,\text{.}
\]

\section{Unified Framework: Rectangular $\ell$-Matchings}

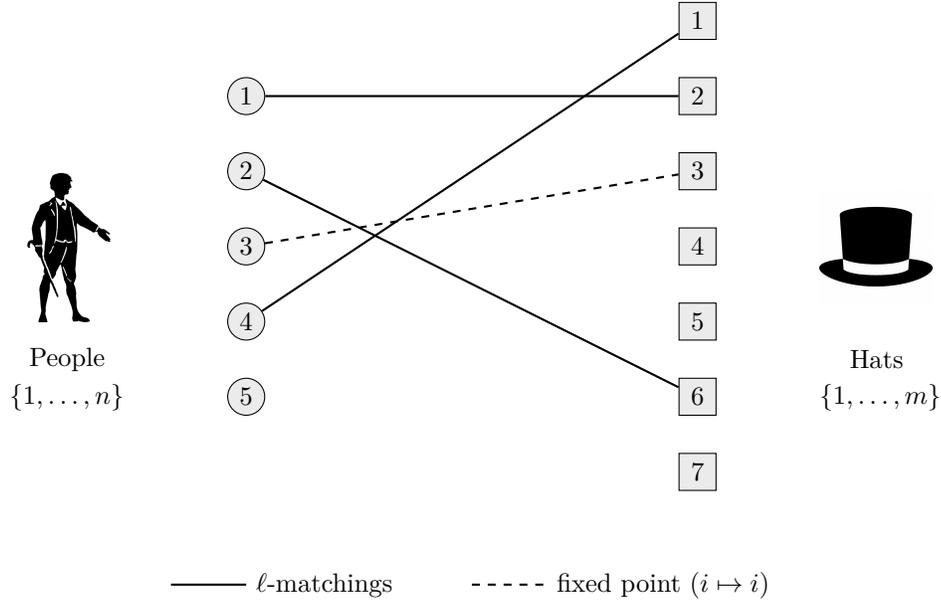
\begin{figure}
\centering
\begin{tikzpicture}[
  >=stealth,
  person/.style={circle, draw, fill=black!8, inner sep=1.5pt, minimum size=14pt},
  hat/.style={rectangle, draw, fill=black!8, inner sep=1.5pt, minimum size=14pt},
  match/.style={thick},
  fixed/.style={thick, dashed}
]

% ====== Parameters to adjust ======
\def\n{5}                         % number of people
\def\m{7}                         % number of hats
\def\assign{{1/2},{2/6},{4/1}, {3/3}}  % list of \ell matchings i/j (here \ell=3)
\def\ptfixe{{1/1}, {2/2}, {3/3}, {4/4}, {5/5}}  % list of fixed points i/i (optional)
\def\xL{0}  % x for people
\def\xR{6}  % x for hats

% ====== Nodes: people (left) ======
\foreach \i in {1,...,\n}{
  \node[person] (p\i) at (\xL,-\i-1) {$\i$};
}
% ====== Nodes: hats (right) ======
\foreach \j in {1,...,\m}{
  \node[hat] (h\j) at (\xR,-\j) {$\j$};
}

% ====== Edges (assignments) ======
\foreach \i/\j in \assign {
  \ifnum\i=\j
    \draw[fixed] (p\i) -- (h\j);   % fixed point: i -> i
  \else
    \draw[match] (p\i) -- (h\j);   % regular assignment
  \fi
}

% ====== Column titles ======
\node[left=17mm] at (\xL+.05 ,-4) {\includegraphics[height=20mm]{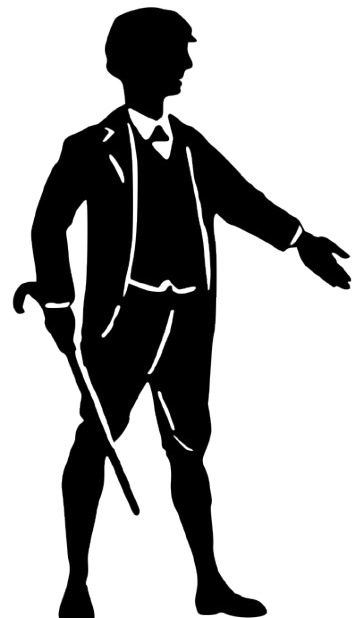}};
\node[left=15mm] at (\xL-.25,-5.5) {People};
\node[left=15mm] at (\xL,-6) {$\{1,\dots, n\}$};

\node[right=16mm] at (\xR-.1 ,-4) {\includegraphics[height=15mm]{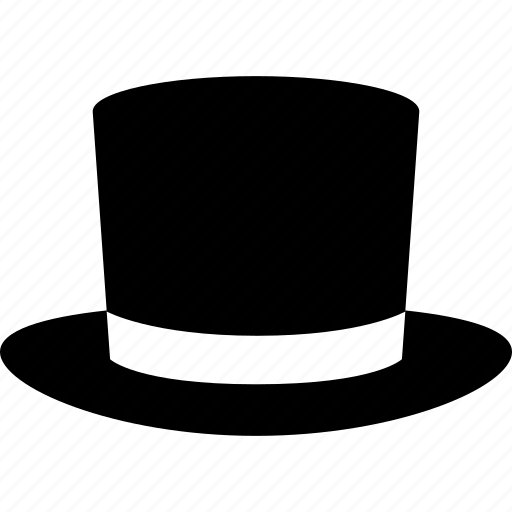}};
\node[right=15mm] at (\xR+.4,-5.5) {Hats};
\node[right=15mm] at (\xR,-6) {$\{1,\dots, m\}$};

% ====== Legend ======
\begin{scope}[shift={(0,-\m-1.5)}]
  \draw[match] (-1,0) -- +(1,0) node[right]{$\ell$-matchings};
  \draw[fixed] (3,0) -- +(1,0) node[right]{fixed point $(i\mapsto i)$};
\end{scope}

\end{tikzpicture}
\caption{Illustration of the general case: \(n\) people (left), \(m\) hats (right), and \(\ell\) matchings (here: \(n=5\), \(m=7\), \(\ell=3\)).}
\label{fig:graphe}
\end{figure}

We can naturally combine the two generalizations and consider the general case with \(n\) people, \(m\) hats (\(m>n\)), and only \(\ell\) of them matched (\(\ell\le n\)).  
A \emph{rectangular \(\ell\)-matching} is an injection \(\sigma:I\hookrightarrow \{1,\dots, m\}\) where \(I\subseteq\{1,\dots, n\}\) and \(|I|=\ell\).  
We denote the set
\[
\grapharrangementsset{n}{m}{\ell}
=\bigl\{\,\sigma:I\hookrightarrow\{1,\dots, m\}\ \text{injection}\ :\ I\subseteq\{1,\dots, n\},\ |I|=\ell\,\bigr\}\,\text{,}
\]
whose cardinality is
\[
\lvert \grapharrangementsset{n}{m}{\ell}\rvert=\binom{n}{\ell}\binom{m}{\ell}\,\ell!\,\text{.}
\]

We define the \emph{partial rectangular derangements} (no fixed point) as
\[
\rectpartialderangements{n}{m}{\ell}
=\bigl\{\,\sigma\in\grapharrangementsset{n}{m}{\ell}:\ \forall i\in I,\ \sigma(i)\neq i\,\bigr\}\,\text{.}
\]
By inclusion–exclusion (Poincaré sieve), we obtain
\[
\bigl\lvert \rectpartialderangements{n}{m}{\ell}\bigr\rvert
=\sum_{j=0}^{\ell}(-1)^j\binom{n}{j}\,\binom{n-j}{\ell-j}\,\binom{m-j}{\ell-j}\,(\ell-j)!\,\text{.}
\]

The probability that a uniform rectangular \(\ell\)-matching has no fixed point is therefore
\[
\Prectpartial{n}{m}{\ell}
=\frac{\bigl\lvert \rectpartialderangements{n}{m}{\ell}\bigr\rvert}
       {\bigl\lvert \grapharrangementsset{n}{m}{\ell}\bigr\rvert}
=\frac{\displaystyle \sum_{j=0}^{\ell}(-1)^j\binom{n}{j}\,\binom{n-j}{\ell-j}\,\binom{m-j}{\ell-j}\,(\ell-j)!}
       {\displaystyle \binom{n}{\ell}\binom{m}{\ell}\,\ell!}\,\text{.}
\]

We then recover the previous special cases:
\begin{itemize}
  \item If \(m=n\), we recover the formula for \(\ell\)-matchings on \(\{1,\dots, n\}\times\{1,\dots, n\}\):  
  \(\binom{n-j}{\ell-j}\binom{m-j}{\ell-j} =  \binom{n-j}{\ell-j}^2\)\,\text{.}
  \item If \(\ell=n\), we recover the case of injections \(\{1,\dots, n\}\hookrightarrow\{1,\dots, m\}\):  
  \(\binom{n-j}{\ell-j}=\binom{n-j}{n-j}=1\) and \(\binom{m-j}{\ell-j}(\ell-j)! = \frac{(m-j)!}{(m-n)!}\,\text{.}\)
\end{itemize}

% \section*{Conclusion}

% Dans cette note, nous proposons d’abord un bref panorama historique — du \emph{Jeu du Treize} de Montmort au \emph{Hat-Check Problem} de Scoville — retraçant l’étude, au fil des siècles, du problème combinatoire des \emph{rencontres}. 
% Nous mettons ensuite en regard le problème initial et deux de ses généralisations naturelles, faisant intervenir des objets de complexité croissante : les permutations (dérangements), les injections rectangulaires (dérangements rectangulaires) et les \(\ell\)-appariements (dérangements partiels). 
% Ces problèmes de dénombrement se traitent de manière unifiée par le \emph{crible de Poincaré} (principe d’inclusion–exclusion), qui conduit à des formules fermées et éclaire la loi du nombre de points fixes (avec une limite de type Poisson).

% Ces généralisations apparaissent dans des contextes variés — par exemple :
% % combinatoire et probabilités (approximation de Poisson, méthode de Chen–Stein),
% théorie des graphes (familles intersectantes, graphes de dérangements, polynômes de tours),
% % algorithmique/crypto (protocoles « Secret Santa », tirages sans auto-affectation),
% % statistique bayésienne (appariements bipartites pour le \emph{record linkage}),
% théorie du codage (codes de permutations sous distance de Hamming).

% À notre connaissance, ces trois cadres sont rarement présentés dans un même système de notations, et plus rarement encore combinés explicitement ; l’objectif de cette note est précisément d’en offrir une présentation unifiée et concise.

\bibliography{sample.bib}

\end{document}